\documentclass[12pt,a4paper]{amsart}
\usepackage{amssymb}
\usepackage{ifthen}
\usepackage{graphicx}
\newtheorem{thm}{Theorem}
\newtheorem{cor}{Corollary}
\newtheorem{lem}{Lemma}

\newtheorem{rem}{Remark}

\newtheorem{conj}{Conjecture}
\theoremstyle{definition}
\newtheorem{defn}[equation]{Definition}%[section]
\newtheorem{example}[equation]{Example}%[section]
\newtheorem{prob}[equation]{Problem}

\newcounter {own}
\def\theown {\thesection       .\arabic{own}}

\newenvironment{pf}[1][]{%
 \vskip 3mm
 \noindent
 \ifthenelse{\equal{#1}{}}%
  {{\slshape Proof. }}%
  {{\slshape #1.} }%
 }%
{\qed\bigskip}

\newcounter{alphabet}
\newcounter{tmp}

\newcommand{\A}{{\mathcal A}}

\newcommand{\U}{{\mathcal U}}

\newcommand{\IC}{{\mathbb C}}
\newcommand{\ID}{{\mathbb D}}

\newcommand{\D}{{\mathbb D}}
\def\be{\begin{equation}}
\def\ee{\end{equation}}

\newcommand{\bee}{\begin{enumerate}}
\newcommand{\eee}{\end{enumerate}}

\newcommand{\blem}{\begin{lem}}
\newcommand{\elem}{\end{lem}}
\newcommand{\bthm}{\begin{thm}}
\newcommand{\ethm}{\end{thm}}
\newcommand{\bcor}{\begin{cor}}
\newcommand{\ecor}{\end{cor}}
\newcommand{\beg}{\begin{example}}
\newcommand{\eeg}{\end{example}}
\newcommand{\begs}{\begin{examples}}
\newcommand{\eegs}{\end{examples}}
\newcommand{\bdefe}{\begin{defin}}
\newcommand{\edefe}{\end{defin}}
\newcommand{\bprob}{\begin{prob}}
\newcommand{\eprob}{\end{prob}}
\newcommand{\bei}{\begin{itemize}}
\newcommand{\eei}{\end{itemize}}

\newcommand{\bcon}{\begin{conj}}
\newcommand{\econ}{\end{conj}}
\newcommand{\bcons}{\begin{conjs}}
\newcommand{\econs}{\end{conjs}}
\newcommand{\bprop}{\begin{propo}}
\newcommand{\eprop}{\end{propo}}
\newcommand{\br}{\begin{rem}}
\newcommand{\er}{\end{rem}}
\newcommand{\brs}{\begin{rems}}
\newcommand{\ers}{\end{rems}}
\newcommand{\bo}{\begin{obser}}
\newcommand{\eo}{\end{obser}}
\newcommand{\bos}{\begin{obsers}}
\newcommand{\eos}{\end{obsers}}
\newcommand{\bpf}{\begin{pf}}
\newcommand{\epf}{\end{pf}}
\newcommand{\ba}{\begin{array}}
\newcommand{\ea}{\end{array}}
\newcommand{\beq}{\begin{eqnarray}}
\newcommand{\beqq}{\begin{eqnarray*}}
\newcommand{\eeq}{\end{eqnarray}}
\newcommand{\eeqq}{\end{eqnarray*}}

\begin{document}
\bibliographystyle{amsplain}

\title[Some properties of the class $\mathcal{U}$]{Some properties of the class $\mathcal{U}$}

\author[M. Obradovi\'{c}]{Milutin Obradovi\'{c}}
\address{Department of Mathematics,
Faculty of Civil Engineering, University of Belgrade,
Bulevar Kralja Aleksandra 73, 11000, Belgrade, Serbia}
\email{obrad@grf.bg.ac.rs}

\author[N. Tuneski]{Nikola Tuneski}
\address{Department of Mathematics and Informatics, Faculty of Mechanical Engineering, Ss. Cyril and Methodius
University in Skopje, Karpo\v{s} II b.b., 1000 Skopje, Republic of Macedonia.}
\email{nikola.tuneski@mf.edu.mk}

%\author{M. Obradovi\'{c}}
%\address{M. Obradovi\'{c},
%Department of Mathematics,
%Faculty of Civil Engineering, University of Belgrade,
%Bulevar Kralja Aleksandra 73, 11000
%\email{obrad@grf.bg.ac.rs}

%\author{S. Ponnusamy${}^{~\mathbf{*}}$}
%\address{S. Ponnusamy, Department of Mathematics,
%Indian Institute of Technology Madras, Chennai--600 036, India.}
%\email{samy@iitm.ac.in}
%\author{N. Tuneski}
%\address{N. Tuneski, St. Cyril and Methodius University, Faculty of Mechanical
%Engineering, Karpo\v s II b.b., 1000 Skopje, R. Macedonia.}
%\email{nikolat@mf.edu.mk}

\subjclass[2000]{30C45, 30C50, 30C55}
\keywords{analytic, class $\mathcal{U}$, starlike, $\alpha$-convex, Hankel determinant.}
%\date{\today  %June. 30, 09
%;  File: }

%\begin{abstract}

%\end{abstract}

%\thanks{The work of the first author was supported by MNZZS Grant, No. ON174017, Serbia. The research of the second
%author was supported by National Board for Higher
%Mathematics, India.}

\maketitle
\pagestyle{myheadings}
%\markboth{M. Obradovi\'{c} and  S. Ponnusamy }{}
%\cc
%\section{Introduction}

\begin{abstract}
In this paper we study the class $\mathcal{U}$ of functions that are analytic in the open unit disk $\D=\{z:|z|<1\}$, normalized such that $f(0)=f'(0)-1=0$ and satisfy
\[\left|\left [\frac{z}{f(z)} \right]^{2}f'(z)-1 \right|<1\quad\quad (z\in\D).\]
For functions in the class $\mathcal{U}$ we give sharp estimate of the second ant the third Hankel determinant, its relationship with the class of $\alpha$-convex functions, as well as certain starlike properties.
\end{abstract}

\bigskip

%As in the previous theorem, let consider first the case when $|c_{1}|^{2}\leq \frac{1}{3}.$ Then we have
%\beqq
%\begin{split}
%|\gamma_{2}|&=\left| c_{2}+\frac{1}{2}c_{1}^{2}\right|\\
%&= \frac{1}{2}\left|(2c_{2}-c_{1}^{2})+2c_{1}^{2}\right|\\
%&\leq \frac{1}{2}|2c_{2}-c_{1}^{2}|+|c_{1}|^{2}\\
%&\leq \frac{1}{2}+\frac{1}{3}=\frac{5}{6}
%\end{split}
%\eeqq
%(where we used inequality in \eqref{eq 95}).
%
%\medskip
%
%In case when $\frac{1}{3}\leq |c_{1}|^{2}\leq 1$, we have
%\beqq
%\begin{split}
%|\gamma_{2}|&=\left| c_{2}+\frac{1}{2}c_{1}^{2}\right|\\
%&\leq |c_{2}|+\frac{1}{2}|c_{1}|^{2}\\
%&\leq 1-|c_{1}|^{2}+\frac{1}{2}|c_{1}|^{2}=1-\frac{1}{2}|c_{1}|^{2}\\
%&\leq1-\frac{1}{2}\cdot\frac{1}{3}=\frac{5}{6}.
%\end{split}
%\eeqq

\section{Introduction}

Let ${\mathcal A}$ denote the family of all analytic functions
in the unit disk $\ID := \{ z\in \IC:\, |z| < 1 \}$ and  satisfying the normalization
$f(0)=0= f'(0)-1$. Let $\mathcal{S}^{\star}$ and $\mathcal{K}$ denote the subclasses of
${\mathcal A}$ which are starlike and convex in $\ID$, respectively, i.e.,
\[ \mathcal{S}^{\star} = \left\{ f\in\mathcal{A}: \operatorname{Re}\left[\frac{zf'(z)}{f(z)}\right] >0,\, z\in\ID \right\} \]
and
\[ \mathcal{K} = \left\{ f\in\mathcal{A}: \operatorname{Re}\left[1+\frac{zf''(z)}{f'(z)}\right] >0,\, z\in\ID \right\}.\]

\medskip

Geometrical characterisation of convexity is the usual one, while for the starlikeness we have that
$f\in \mathcal{S}^\star$, if, and only if, $f(\D)$ is a starlike region, i.e.,
\[z\in f(\D)\quad\quad \Rightarrow  \quad\quad tz\in f(\D) \mbox{ for all } t\in[0,1] .\]

\medskip

The linear combination of the expressions involved in the analytical representations of starlikeness and convexity brings us to
the classes of $\alpha$-convex functions introduced in 1969 by Mocanu (\cite{Mocanu-1969}) and consisting of functions   $f\in \mathcal{A}$ such that
 \be\label{eq 6}
{\operatorname{Re}}\left\{(1-\alpha)\frac{zf'(z)}{f(z)}+ \alpha \left[1+ \frac{z f''(z)}{f'(z)}\right]\right\}>0, \quad\quad (z\in \ID),
\ee
where $\frac{f(z)f'(z)}{z}\neq 0$ for  $z \in \ID $ and $\alpha\in\mathbb{R}$.
Those classes he denoted by $\mathcal{M_{\alpha}}$.

\medskip

Further, let $\mathcal{U} $ denote the set of all  $f\in {\mathcal A}$ satisfying
the condition
$$\left |\operatorname U_f(z)\right | < 1 \quad\quad (z\in \ID),$$
where the operator $\operatorname U_f$ is defined by
\[\operatorname U_f(z):= \left [\frac{z}{f(z)} \right]^{2}f'(z)-1.\]

\medskip

All this classes consist of univalent functions and more details on them can be found in \cite{Good-1,TTV}.

\medskip

The class of starlike functions is very large and in the theory of univalent functions it is significant if a class doesn't entirely lie inside $\mathcal{S}^{\star}$. One such case is the class of functions with bounded turning consisting of functions $f$ from $\mathcal{A}$ that satisfy $\operatorname{Re}f'(z)>0$ for all $z\in \D$. Another example is the class $\mathcal{U}$ defined above and first treated in \cite{Japonica_1996} (see also \cite{obpon-1,obpon-2,TTV}). Namely, the function $-\ln(1-z)$ is convex, thus starlike, but not in $\mathcal{U}$ because $\operatorname U_f(0.99)=3.621\ldots>1$, while  the function $f$ defined by $\frac{z}{f(z)}=1-\frac32 z + \frac12z^{3}=(1- z)^{2} \left(1+\frac{z}{2}\right)$ is in $\mathcal{U}$ and such that $\frac{zf'(z)}{f(z)} = -\frac{2 \left(z^2+z+1\right)}{z^2+z-2} = -\frac{1}{5}+\frac{3 i}{5}$ for $z=i$. This rear property is the main reason why the class $\mathcal{U}$ attracts huge attention in the past decades.

\medskip

In this paper we give sharp estimates of the second and the third Hankel determinant over the class $\mathcal{U}$ and study its relation with the class of $\alpha$-convex and starlike functions.

\medskip

\section{Main results}

\medskip

In the first theorem we give the sharp estimates of the Hankel determinants of second and third order
for the class $\mathcal{U}$. We first give the definition of the Hankel determinant, whose elements are the coefficients of a function $f\in \A$.

\medskip

\begin{defn}
Let $f\in \A$. Then the $qth$ Hankel determinant of $f$ is defined for $q\geq 1$, and
$n\geq 1$ by
\[
        H_{q}(n) = \left |
        \begin{array}{cccc}
        a_{n} & a_{n+1}& \ldots& a_{n+q-1}\\
        a_{n+1}&a_{n+2}& \ldots& a_{n+q}\\
        \vdots&\vdots&~&\vdots \\
        a_{n+q-1}& a_{n+q}&\ldots&a_{n+2q-2}\\
        \end{array}
        \right |.
\]
\end{defn}
Thus, the second and the third Hankel determinants are, respectively,
\be\label{eq 1}
\begin{split}
H_{2}(2)&= a_2a_4-a_{3}^2,\\
H_{3}(1)&= a_{3}(a_2a_4-a_{3}^2)-a_{4}(a_4-a_{2}a_{3})+a_{5}(a_3-a_{2}^2).
\end{split}
\ee

\medskip

\bthm \label{th 1}
Let $f\in \mathcal{U}$ and  $f(z)=z+a_{2}z^{2}+a_{3}z^{3}+\ldots$. Then we have the sharp estimates:
$$|H_{2}(2)|\leq1\quad \mbox{and}\quad  |H_{3}(1)|\leq \frac{1}{4}.$$
\ethm

\medskip

\begin{proof}
In \cite{Japonica_1996} the following characterization for functions $f$ in the class in $\U$ was given:
\be\label{eq 2}
\frac{z}{f(z)} = 1-a_2 z -z\int_0^z \frac{\omega(t)}{t^2}\,dt,
\ee
where function $\omega$ is analytic in $\D$ with $\omega(0)=\omega'(0)=0$ and $|\omega(z)|<1$ for all $z\in\D$.

\medskip

If we put $\omega_1(z)=\int_0^z \frac{\omega(t)}{t^2}\,dt$, then we easily obtain that $|\omega_1(z)|\le|z|<1$ and $|\omega_1'(z)|\le1$ for all $z\in\D$. If $\omega_1(z)=c_1z+c_2z^2+\cdots$, then $\omega_1'(z)=c_1+2c_2z+3c_3z^2+\cdots$ and $|\omega_1'(z)|\le1$, $z\in\D$, gives (see relation (13) in the paper of Prokhorov and Szynal \cite{Prokhorov-1984}):
\be\label{eq 3}
|c_{1}|\leq 1,\quad |2c_2|\leq1-|c_1|^2\quad\mbox{and}\quad |3c_{3}(1-|c_{1}|^{2})+4\overline{c_{1}}c_{2}^{2}|\leq (1-|c_{1}|^{2})^{2}- 4|c_2|^{2}.
\ee
Also, from \eqref{eq 2} we have
\beqq
\begin{split}
 f(z) &= \frac{z}{1-\left(a_2z+c_1z^2+c_2z^3+\cdots\right)} \\
 &= z+a_2z^2+\left(c_1+a_2^2\right)z^3+\left(c_2+2a_2c_1+a_2^3\right)z^4\\
 &\phantom{= }+ \left(c_3+2a_2c_2+c_{1}^{2}+3a_2^2 c_{1}+a_{2}^{4}\right)z^5 \cdots.
\end{split}
\eeqq
From the last relation we have
\be\label{eq 4}
 a_3 = c_1+a_2^2,\quad a_4 = c_2+2a_2c_1+a_2^3,\quad a_{5}=c_3+2a_2c_2+c_{1}^{2}+3a_2^2 c_{1}+a_{2}^{4}.
 \ee
We may suppose that $c_{1}\geq 0$, since from \eqref{eq 4} we have $c_{1}=a_{3}-a_{2}^{2}$ and $a_{3}$ and
$a_{2}^{2}$ have the same turn under rotation. In that sense, from \eqref{eq 3} we obtain
\be\label{eq 5}
0\leq c_{1}\leq 1,\quad |c_2|\leq\frac{1}{2}\left(1-c_1^2\right)\quad \mbox{and}\quad   |c_{3}|\leq\frac{1}{3}\left(1-c_{1}^{2}-\frac{4|c_{2}|^{2}}{1+c_{1}}\right).
\ee
If we use \eqref{eq 1}, \eqref{eq 4} and \eqref{eq 5}, then
\beqq
\begin{split}
 \left|H_{2}(2)\right| =& \left| c_2a_2-c_1^2 \right| \le |c_2|\cdot |a_2| + c_1^2 \le \frac12 \left(1-c_1^2\right) |a_2| +c_1^2 \\
 =& \frac12\cdot|a_2|+\left( 1-\frac12\cdot|a_2| \right)c_1^2 \le1.
 \end{split}
\eeqq

\medskip

The functions $k(z)=\frac{z}{(1-z)^2}$ and $f_{1}(z)=\frac{z}{1-z^2}$ show that the estimate is the best possible.

\medskip

Similarly, after some calculations we also have
\beqq
\begin{split}
\left|H_{3}(1)\right|&= \left|c_{1}c_{3}-c_{2}^{2}\right|\leq c_{1}|c_{3}|+|c_{2}|^{2}\\
&\leq \frac{1}{3}c_{1}\left(1-c_{1}^{2}-\frac{4|c_{2}|^{2}}{1+c_{1}}\right)+|c_{2}|^{2}\\
&=\frac{1}{3}\left(c_{1}-c_{1}^{3}+\frac{3-c_{1}}{1+c_{1}}|c_{2}|^{2}\right)\\
&= \frac{1}{3}\left(c_{1}-c_{1}^{3}+\frac{3-c_{1}}{1+c_{1}}\cdot\frac{1}{4}\left (1-c_{1}^{2}\right)^{2}\right)\\
&= \frac{1}{12}\left(3-2c_{1}^{2}-c_{1}^{4}\right)\leq\frac{3}{12}=\frac{1}{4}.
\end{split}
\eeqq
The function $f_{2}(z)=\frac{z}{1-z^{3}/2}$ shows that the result is the best possible.
\end{proof}

In the rest of the paper be consider some starlikeness problems for the class $\mathcal{U}$ and its connection with the class of $\alpha$-convex functions.

\medskip

First, let recall the classical results about the relation between the starlike functions and $\alpha$-convex functions.

\medskip

\bthm\label{th-A} $ $
\begin{itemize}
  \item[(a)] $\mathcal{M_{\alpha}}\subseteq \mathcal{S}^{\star}$ for every real $\alpha$ $($\cite{Mocanu-1973}$)$;
  \item[(b)] for $0\leq \frac{\beta}{\alpha}\leq1$ we have $\mathcal{M_{\alpha}} \subset \mathcal{M_{\beta}}$  and for $\alpha > 1$,  $\mathcal{M_{\alpha}}\subset\mathcal{M}_{1}=\mathcal{K}$ $($\cite{saka,Mocanu-1973}$)$.
\end{itemize}
\ethm

\medskip

As an anlogue of the above theorem we have

\medskip

\bthm \label{th 2}
For the classes $\mathcal{M_{\alpha}}$ the next results are valid.
\begin{itemize}
  \item[(a)] $\mathcal{M_{\alpha}}\subset\mathcal{U}$ for $\alpha \leq -1$;
  \item[(b)] $\mathcal{M_{\alpha}}$ is not a subset of  $\mathcal{U}$ for any  $0\leq \alpha \leq 1$.
\end{itemize}
\ethm

\begin{proof}$ $
\begin{itemize}
\item[(a)] Let $p(z)=\operatorname U_{f}(z). $
Then $p$ is analytic in $\ID$
and $p(0)=p'(0)=0$. From here we have that $\left [\frac{z}{f(z)} \right ]^{2}f'(z)=p(z)+1$ and, after some calculations that
$$2\frac{zf'(z)}{f(z)}-\left[1+ \frac{z f''(z)}{f'(z)}\right]=1-\frac{zp'(z)}{p(z)+1}.$$
The relation \eqref{eq 6} is equivalent to
\be\label{eq 7}
{\operatorname{Re}}\left\{(1+\alpha)\frac{zf'(z)}{f(z)}-\alpha \left[ 1-\frac{zp'(z)}{p(z)+1}\right]\right\}>0,\,\,z\in \ID.
\ee
We want to prove that $|p(z)|<1$, $z\in \ID$. If not, then according to the Clunie-Jack Lemma (\cite{jack}) there exists a $z_{0}$, $|z_{0}|<1$, such that
$p(z_{0})=e^{i\theta}$ and $z_{0}p'(z_{0})=kp(z_{0})=ke^{i\theta}$, $k \geq 2.$ For such $z_{0}$, from
\eqref{eq 7} we have that
\[
\begin{split}
{\operatorname{Re}}&\left\{(1+\alpha)\frac{z_{0}f'(z_{0})}{f(z_{0})}-\alpha \left[ 1-\frac{ke^{i\theta}}{e^{i\theta}+1}\right]\right\}\\
=& (1+\alpha){\operatorname{Re}}\left[\frac{z_{0}f'(z_{0})}{f(z_{0})}\right]+\alpha \frac{k-2}{2}\leq 0
\end{split}
\]
since $f\in \mathcal{S}^{\star}$ (by Theorem \ref{th-A}) and $\alpha \leq -1$. That is the
contradiction to \eqref{eq 6}. It means that $|p(z)|=|\operatorname U_{f}(z)|<1$, $z\in \ID,$
i.e. $f\in \mathcal {U}.$
\medskip
  \item[(b)] To prove this part, by using Theorem \ref{th-A}(b), it is enough to find a function $ g\in \mathcal{K}$ such that
$g $ not belong to the class $\mathcal {U}$.
Really, the function $g(z)=-\ln(1-z)$ is convex but not in $\mathcal {U}.$
\end{itemize}
\end{proof}

\medskip

\noindent
{\bf Open problem}. It remains an open problem to study the relationship between classes $\mathcal{M_{\alpha}}$ and $\mathcal{U}$ when $-1<\alpha<0$ and $\alpha >1$.

\bigskip

In the next theorem  we  consider starlikeness of the function
\be\label{eq 8}
g(z)=\frac{z/f(z)-1}{-a_{2}}
\ee
where $f\in\mathcal{U}$  and   $a_{2}=\frac{f''(0)}{2}\neq 0$, i.e., its second coefficient doesn't vanish.

\medskip

Namely, we have

\bthm\label{th 3}
Let $f\in\mathcal{U}$. Then, for the function $g$ defined by \eqref{eq 8} we have:
\begin{itemize}
  \item[$(a)$] $|g'(z)-1|<1$ for $|z|<|a_{2}|/2$;
  \item[$(b)$] $g\in \mathcal{S}^{\star}$ in the disc  $|z|<|a_{2}|/2$ and even more
  \[ \left|\frac{zg'(z)}{g(z)}-1\right| <1 \quad\quad (|z|<|a_{2}|/2);\]
  \item[$(c)$] $g\in \mathcal{U}$ in the disc  $|z|<|a_{2}|/2$ if $ 0<|a_{2}|\leq 1$.
\end{itemize}
The results  are best possible.
\ethm

\begin{proof} Let $f\in\mathcal{U}$ with  $a_{2}\neq 0$ . Then, by using \eqref{eq 2}, we have that
$$\frac{z}{f(z)}=1-a_{2}z-z\omega_{1}(z),$$
where $\omega_{1} $ is analytic in $\ID$ such that
$|\omega_{1}(z)|\leq|z| $ and $|\omega'_{1}(z)|\leq1 $. The appropriate function $g$ from \eqref{eq 8} has the form
\[
g(z)=z+\frac{1}{a_{2}}z\omega_{1}(z).
\]
From here $|g'(z)-1|=\frac{1}{|a_2|}|\omega_1(z)+ z\omega_1'(z)|<1 $ for $|z|<|a_{2}|/2$.

\medskip

By using previous representation, we obtain
\beqq
\left|\frac{zg'(z)}{g(z)}-1\right| = \left|\frac{z\omega'_{1}(z)}{a_{2}+\omega_{1}(z)}\right| \le \frac{|z|}{|a_2|-|z|}< 1
\eeqq
if $|z|<|a_{2}|/2$. It means that the function $g$ is starlike in the disk $|z|<|a_{2}|/2$.

\medskip

If we consider function $f_{b}$ defined by
\be\label{eq 10}
\frac{z}{f_{b}(z)}=1+bz+z^{2},\quad \quad 0< b\leq2,
\ee
then $f_{b}\in\mathcal{U}$ and $$ g_{b}(z)=\frac{\frac{z}{f_{b}(z)}-1}{b}=z+\frac{1}{b}z^{2}.$$
For this function we easily have that for $|z|<b/2$:
$$ {\operatorname{Re}}\frac{zg'_{b}(z)}{g_{b}(z)}\geq\frac{1-\frac{2}{b}|z|}{1-\frac{1}{b}|z|}>0.$$
On the other hand side, since
$g'_{b}(-b/2)=0$, the function $g_{b}$ is not univalent in a bigger disc, which implies that our result is  best possible.

\medskip

Also, by using \eqref{eq 8} and the next estimation for the function $\omega_{1}$:
$$|z\omega'_{1}(z)-\omega_{1}(z)|\leq\frac{r^{2}-|\omega_{1}(z)|^{2}}{1-r^{2}},$$ (where $|z|=r$ and
$|\omega_{1}(z)|\leq r$), after some calculation we get
\beqq
\begin{split}
|\mathcal{U}_{g}(z)|&=
\left|\frac{\frac{1}{a_{2}}(z\omega'_{1}(z)-\omega_{1}(z))-\frac{1}{a^{2}_{2}}\omega^{2}_{1}(z)}{\left(1+\frac{1}{a_{2}}\omega_{1}(z)\right)^{2}}\right|\\
&\leq\frac{|a_{2}||z\omega'_{1}(z)-\omega_{1}(z)|+|\omega_{1}(z)|^{2}}{(|a_{2}|-|\omega_{1}(z)|)^{2}}\\
&\leq \frac{|a_{2}|\frac{r^{2}-|\omega_{1}(z)|^{2}}{1-r^{2}}+|\omega_{1}(z)|^{2}}{(|a_{2}|-|\omega_{1}(z)|)^{2}}\\
&=:\frac{1}{1-r^{2}}\varphi(t),
\end{split}
\eeqq
where we put $$\varphi(t)=\frac{(1-r^{2}-|a_{2}|)t^{2}+|a_{2}|r^{2}}{(|a_{2}|-t)^{2}}$$
and $|\omega_{1}(z)|=t$, $0\leq t\leq r$. Since $$\varphi'(t)=\frac{2|a_{2}|}{(|a_{2}|-t)^{3}}\left((1-r^{2}-|a_{2}|)t+r^{2}\right)
=\frac{2|a_{2}|}{(|a_{2}|-t)^{3}}\left((1-|a_{2}|)t+(1-t)r^{2}\right)\geq 0,$$
because $ 0<|a_{2}|\leq 1$ and $0\leq t<1$. It means that the function
 $\varphi$ is an increasing function and that
$$\varphi(t)\leq\varphi(r)=\frac{(1-r^{2})r^{2}}{(|a_{2}|-r)^{2}}.$$

\medskip

Finally we have that
$$|\operatorname U_{g}(z)|\leq\frac{r^{2}}{(|a_{2}|-r)^{2}}<1,$$
since $|z|<|a_{2}|/2$. That is implies the second statement of the theorem.

\medskip

As for sharpness, we can also consider the function  $f_{b}$ defined by \eqref{eq 10} with $0<b\leq 1$.
For $|z|<\frac{b}{2}$ we have
$$\left|\operatorname U_{g_{b}}(z)\right|\leq\frac{\frac{1}{b^{2}}|z|^{2}}{\left(1-\frac{1}{b}|z|\right)^{2}}<1,$$
which implies that $g_{b}$
belongs to the class $\mathcal{U}$ in the disc  $|z|<b/2$.
\end{proof}

\medskip

We believe that part (b) of the previous theorem is valid for all $0<|a_{2}|\leq 2$.
In that sense we have the next

\begin{conj}Let $f\in\mathcal{U}$. Then the function $g$ defined by the expression \eqref{eq 8}
belongs to the class  $\mathcal{U}$ in the disc  $|z|<|a_{2}|/2$. The result is the best possible.
\end{conj}

\medskip

\end{document}